\documentclass[10pt,twoside]{article}
\usepackage{Latex-document-web}
\almvol{00}
\almttone{}%Frontiers of Science Awards}
\almtttwo{}%for Math/TCIS/Phys}
\firstpage{1}
\usepackage[english]{babel}
\usepackage[utf8]{inputenc}

\usepackage{multicol}

\usepackage{tabu}

\usepackage{amsmath}
\usepackage{amssymb}

\usepackage{tikz}
\usetikzlibrary{shapes}
\usetikzlibrary{calc}
\usepackage{pgfplots}
\usepgfplotslibrary{fillbetween}

\numberwithin{equation}{section}

\newtheorem{thm}{Theorem}[section]

\newtheorem{con}[thm]{Conjecture}

\newtheorem{prob}[thm]{Problem}

\begin{document}

\markboth{\hfill{\rm David Conlon} \hfill}{\hfill {\rm Combinatorial theorems relative to sparse sets \hfill}}

\title{Combinatorial theorems\\ relative to sparse sets}

\author{David Conlon}

\begin{abstract}
A key theme in modern extremal combinatorics is the study of classical combinatorial theorems relative to sparse subsets of their natural settings. Here we describe some of the recent progress in this area and state a number of problems that remain open and pressing. 
\end{abstract}

\maketitle

\setcounter{tocdepth}{1}
\tableofcontents

\section{Introduction}

In 1986, reproving a result of Folkman~\cite{F70}, Frankl and R\"odl~\cite{FR86} showed that there are $K_4$-free graphs with the property that every two-colouring of their edges contains a monochromatic $K_3$. In the same paper, they also proved the existence of $K_4$-free graphs $G$ with the property that every subgraph of $G$ containing $(\frac 12 + \epsilon) e(G)$ edges contains a $K_3$, where, because there are always bipartite subgraphs with at least $\frac 12 e(G)$ edges, $\frac 12$ is best possible. Underlying these results are proofs that sparse random graphs satisfy versions of the triangle cases of both Ramsey's theorem and Tur\'an's theorem, two of the basic results in extremal combinatorics. From this initial seed has grown an entire area studying when combinatorial theorems hold relative to sparse sets. It is this area, and its many facets, that we discuss here. 

In the random case, long after an important early success by R\"odl and Ruci\'nski~\cite{RR93, RR95} regarding Ramsey properties of random graphs, significant progress in our understanding was made roughly a decade ago when Conlon and Gowers~\cite{CG16} and, independently, Schacht~\cite{S16} (see also~\cite{FRS10}) developed general techniques for `transferring' combinatorial theorems from their original settings to the random setting. This was soon followed by the development of the hypergraph container method by Balogh, Morris and Samotij~\cite{BMS15} and Saxton and Thomason~\cite{ST15}, allowing another approach to these problems and much more besides. Alongside this, stemming in part from Green and Tao's seminal work~\cite{GT08} on arithmetic progressions in the primes, there has been considerable work on combinatorial theorems relative to pseudorandom sets. More recently, there has been work on combinatorial theorems relative to extremal sets such as $C_4$-free graphs and Sidon sets, with surprising consequences in extremal and additive combinatorics. We will discuss each of these topics, random, pseudorandom and extremal, in turn.

Since there are already several surveys discussing these and related topics~\cite{BMS18, B17, C14, CFZ14b, GS05, RS13}, we will not attempt to be comprehensive here, focusing instead on some of the more recent developments and posing or reiterating a number of important problems that remain open. We begin by looking at the random case. 

\section{Combinatorial theorems relative to random sets}

Given a finite set $X$ and $0 < p \le 1$, let $X_p$ be the random subset formed by placing each element of $X$ in $X_p$ independently with probability $p$. A much-studied special case of this construction is the \emph{binomial random graph} $G_{n,p}$ where each edge of $K_n$ is chosen to be in $G_{n,p}$ independently with probability $p$. The theme of this section will be the question of determining for which probabilities $p$ a given combinatorial theorem about a finite set $X$ remains true with high probability relative to the random subset $X_p$.

Despite the fact that it is thirty years since it was proved, the exemplar of such a theorem is still the random Ramsey theorem of R\"odl and Ruci\'nski~\cite{RR93, RR95}, which determines the threshold for a random graph to be Ramsey with respect to a fixed graph $H$. 
Formally, given a graph $H$ and a natural number $r \ge 2$, we say that a graph $G$ is \emph{$(H,r)$-Ramsey} if every $r$-colouring of the edges of $G$ contains a monochromatic copy of $H$. Ramsey's theorem itself is the statement that $K_n$ is $(H,r)$-Ramsey provided $n$ is sufficiently large in terms of $H$ and $r$. The R\"odl--Ruci\'nski theorem is then as follows, where, here and throughout this survey, we write $v(H)$ and $e(H)$ for the number of vertices and edges of a graph $H$.

\begin{thm} \label{thm:RR}
For any graph $H$ that is not a forest consisting of stars and paths of length $3$ and any integer~$r \geq 2$, there exist positive constants $c$ and $C$ such that 
\[
\lim_{n \rightarrow \infty} \mathbb{P} [G_{n,p} \mbox{ is $(H,r)$-Ramsey}] =
\begin{cases}
0 & \text{if $p < c n^{-1/m_2(H)}$}, \\
1 & \text{if $p > C n^{-1/m_2(H)}$},
\end{cases}
\]
where 
\[m_2(H) = \max\left\{\frac{e(H') - 1}{v(H') - 2}: H' \subseteq H \mbox{ and } v(H') \geq 3\right\}.\]
\end{thm}

This theorem is really two theorems in one, each of which was proved in a separate paper. The first~\cite{RR93} is the $0$-statement, stating that the probability $G_{n,p}$ is $(H,r)$-Ramsey is asymptotically $0$ when $p < c n^{-1/m_2(H)}$, while the other~\cite{RR95} is the $1$-statement, saying that for $p > C n^{-1/m_2(H)}$ the probability is asymptotically $1$. The original proofs of both of these results were quite difficult, though both have been considerably simplified in recent years~\cite{NS16}. 

An important recent development, resolving a longstanding conjecture of Kohayakawa and Kreuter~\cite{KK97}, extends the R\"odl--Ruci\'nski theorem to the asymmetric setting. We say that a graph is \emph{$(H_1, \dots, H_r)$-Ramsey} if every $r$-colouring of the edges of $G$ contains a copy of $H_i$ in colour $i$ for some $i$. The final result, a combination of work by Mousset, Nenadov and Samotij~\cite{MNS20}, Bowtell, Hancock and Hyde~\cite{BHH25+}, Kuperwasser, Samotij and Wigderson~\cite{KSW25} and Christoph, Martinsson, Steiner and Wigderson~\cite{CMSW25}, is as follows.

\begin{thm}
For any graphs $H_1, \dots, H_r$ with $m_2(H_1) \ge m_2(H_2) \ge \dots \ge m_2(H_r)$ and $m_2(H_2) > 1$, there exist positive constants $c$ and $C$ such that 
\[
\lim_{n \rightarrow \infty} \mathbb{P} [G_{n,p} \mbox{ is $(H_1,\dots, H_r)$-Ramsey}] =
\begin{cases}
0 & \text{if $p < c n^{-1/m_2(H_1, H_2)}$}, \\
1 & \text{if $p > C n^{-1/m_2(H_1,H_2)}$},
\end{cases}
\]
where
\[m_2(H_1, H_2) = \max\left\{\frac{e(H')}{v(H') - 2 + 1/m_2(H_2)}: H' \subseteq H_1 \mbox{ and } v(H') \geq 2\right\}.\]
\end{thm}

Despite these successes, the situation for hypergraphs remains somewhat mysterious. Though an analogue of the $1$-statement in the R\"odl--Ruci\'nski theorem is known~\cite{CG16, FRS10}, there are more situations than in the graph case where this is not the sharp bound~\cite{GNPSST17}. While it might be very difficult, or even impossible, to determine the threshold for the Ramsey property for hypergraphs in full generality, our knowledge of the graph case is now sufficiently advanced that it seems a reasonable goal to study the $3$-uniform case in more detail.

\begin{prob}
Determine the threshold for the $(H,r)$-Ramsey property for all $3$-uniform hypergraphs $H$ and all $r \ge 2$.
\end{prob}

Another interesting open problem, though perhaps lying further away, is to show that the threshold for the Ramsey property, even for $H = K_3$ and $r = 2$, has a sharp threshold of the form $C/\sqrt{n}$. The fact that the threshold is sharp, in the sense that there exists some function $p(n)$, not necessarily of the form $C/\sqrt{n}$, such that the threshold concentrates around this function, for $H = K_3$ and $r = 2$ is a theorem of Friedgut, R\"odl, Ruci\'nski and Tetali~\cite{FRRT06}, making fundamental use of a seminal result of Friedgut~\cite{F99} characterising sharp thresholds for monotone properties. This result has been considerably extended since by Schacht and Schulenburg~\cite{SS18} and by Friedgut, Kuperwasser, Samotij and Schacht~\cite{FKSS25+}. In particular, the latter paper shows that the $(H,r)$-Ramsey property has a sharp threshold for $H$ any clique or cycle and $r$ any number of colours. However, the problem of determining the true form of the threshold remains wide open (and may be very difficult).

\begin{prob}
Show that the $(K_3, 2)$-Ramsey property has a sharp threshold at $C/\sqrt{n}$ for some $C > 0$.
\end{prob}

An analogue of Theorem~\ref{thm:RR} for the Tur\'an property was conjectured by Haxell, Kohayakawa and \L uczak~\cite{HKL95, HKL96} in the mid 1990s.  
Given a graph $H$ and $\epsilon > 0$, we say that a graph $G$ is \emph{$(H,\epsilon)$-Tur\'an} if every subgraph of $G$ with at least $\left(1 - \frac{1}{\chi(H)-1} + \epsilon\right) e(G)$ edges contains a copy of $H$. By the Erd\H{o}s--Stone--Simonovits theorem (whose special case with $H = K_t$ is essentially Tur\'an's theorem), we know that provided $n$ is sufficiently large in terms of $H$ and $\epsilon$ the complete graph $K_n$ is $(H,\epsilon)$-Tur\'an. The conjecture of Haxell, Kohayakawa and \L uczak, now a theorem proved independently by Conlon and Gowers~\cite{CG16} and by Schacht~\cite{S16}, states that the Tur\'an property holds down to essentially the same threshold as the Ramsey property. There is one caveat here, which is that the result of Conlon and Gowers applied only to \emph{strictly $2$-balanced} graphs, those $H$ such that $m_2(H) > m_2(H')$ for all proper subgraphs $H'$ of $H$. However, almost all graphs satisfy this requirement, including all cliques and cycles.

\begin{thm} \label{thm:randomTuran}
For any graph $H$ and any $\epsilon > 0$, there exist positive
constants $c$ and $C$ such that
\[
\lim_{n \rightarrow \infty} \mathbb{P} [G_{n,p} \mbox{ is $(H,\epsilon)$-Tur{\'a}n}] =
\begin{cases}
0 & \text{if $p < c n^{-1/m_2(H)}$}, \\
1 & \text{if $p > C n^{-1/m_2(H)}$}.
\end{cases}
\]
\end{thm}

We note that the $0$-statement here is surprisingly straightforward. Indeed, if the number of copies of $H$ is significantly smaller than the number of edges, we can remove all copies of $H$ by deleting one edge from each copy. Therefore, if $p^{e(H)} n^{v(H)} \ll p n^2$, that is, $p \ll n^{-(v(H)-2)/(e(H)-1)}$, the $(H, \epsilon)$-Tur\'an property cannot hold. Since the same argument applies for any subgraph $H'$ of $H$, it is easy to see that for $p \ll n^{-1/m_2(H)}$ the random graph $G_{n,p}$ cannot be $(H, \epsilon)$-Tur\'an. Hence, all of the difficulty lies in proving the $1$-statement.

The results of~\cite{CG16, S16} (see also~\cite{CGSS14, FRS10, S14}) prove much more than Theorem~\ref{thm:randomTuran}, allowing one to `transfer' many combinatorial theorems about bounded-size objects from their natural setting to a sparse random subset of that setting. Another representative example is the following random version of Szemer\'edi's theorem~\cite{Sz75} on arithmetic progressions in dense sets of integers. To state the result, given an integer $k \ge 3$ and $\delta > 0$, we say that a subset $I$ of the integers is \emph{$(k,\delta)$-Szemer\'edi} if any subset of $I$ with at least $\delta |I|$ elements contains an arithmetic progression of length $k$. Szemer\'edi's theorem says that for $n$ sufficiently large in terms of $k$ and $\delta$ the set $[n] := \{1, 2, \dots, n\}$ is $(k, \delta)$-Szemer\'edi, while a striking corollary of Green and Tao's work on arithmetic progressions in the primes~\cite{GT08} says that for $n$ sufficiently large in terms of $k$ and $\delta$ the set of primes up to $n$ is $(k, \delta)$-Szemer\'edi. 

\begin{thm} \label{thm:Szem}
For any integer $k \geq 3$ and $\delta > 0$, there exist positive constants $c$ and $C$ such that 
\[
\lim_{n \rightarrow \infty} \mathbb{P} [[n]_{p} \mbox{ is $(k,\delta)$-Szemer\'edi}] =
\begin{cases}
0 & \text{if $p < c n^{-1/(k-1)}$}, \\
1 & \text{if $p > C n^{-1/(k-1)}$}.
\end{cases}
\]
\end{thm}

Another example is the graph removal lemma, a key corollary of the celebrated regularity lemma~\cite{Sz78} due to Ruzsa and Szemer\'edi~\cite{RSz78} stating that any $n$-vertex graph with $o(n^{v(H)})$ copies of a fixed graph $H$ can be made $H$-free by removing $o(n^2)$ edges. This is a fundamental result in extremal combinatorics, implying, amongst other things, Roth's theorem~\cite{R53}, the $3$-term progression case of Szemer\'edi's theorem. The fact that it also holds in the sparse context was first proved for strictly $2$-balanced graphs by Conlon and Gowers~\cite{CG16} and then extended to all graphs by Conlon, Gowers, Samotij and Schacht~\cite{CGSS14}.

\begin{thm} \label{thm:removal-Gnp}
For any graph $H$ and any $\epsilon > 0$, there exist positive constants $\delta$ and $C$ such that if $p \geq Cn^{-1/m_2(H)}$, then the random graph $G_{n,p}$ a.a.s.\ has the property that every subgraph of $G_{n,p}$ which contains at most $\delta p^{e(H)} n^{v(H)}$ copies of $H$ may be made $H$-free by removing at most $\epsilon p n^2$ edges.
\end{thm}

It is worth remarking that hypergraph analogues of these results also hold. For instance, even though we do not know the Tur\'an density of most fixed hypergraphs, the results of~\cite{CG16, S16} imply that whatever the Tur\'an density of a given $k$-uniform hypergraph $H$ is, it is also the Tur\'an density of $H$ relative to the binomial random hypergraph $G_{n,p}^{(k)}$ down to $p$ roughly $n^{-1/m_k(H)}$, where 
\[m_k(H) = \max\left\{\frac{e(H') - 1}{v(H') - k}: H' \subseteq H \mbox{ and } v(H') \geq k+1\right\}.\]
Similarly, for the hypergraph removal lemma, an important result of Gowers~\cite{G07} and Nagle, R\"odl, Schacht and Skokan~\cite{NRS06, RS04} known to imply Szemer\'edi's theorem in full generality, there is a sparse analogue which holds down to the same threshold. We will say more about sparse hypergraph removal lemmas in the pseudorandom context below.

The approaches taken by Conlon and Gowers~\cite{CG16} and by Schacht~\cite{S16} are quite different and each has its own strengths and weaknesses. Remarkably, not long after their results were made public, a third approach, the hypergraph container method, was discovered independently by two groups of authors, Balogh, Morris and Samotij~\cite{BMS15} and Saxton and Thomason~\cite{ST15}. This method has been remarkably influential, with myriad applications across extremal and probabilistic combinatorics and beyond. We will not attempt to do justice to this topic here, referring the reader instead to the survey paper~\cite{BMS18} and its references. However, it is worth noting that, unlike the other methods (though see~\cite{CGSS14} for a variant), this technique allowed a full resolution of the K\L R conjecture~\cite{KLR97}, a technical statement which enables the application of regularity methods in the sparse random setting and which was seen, for many years, as the most likely route towards proving statements like Theorem~\ref{thm:randomTuran}. We refer the reader also to~\cite{N22} for a surprisingly simple recent proof of the K\L R conjecture by Nenadov.

Though there are several problems where additional ideas are needed (see, for example,~\cite{KS25}), one might reasonably claim that these results effectively resolve the problem of transferring combinatorial theorems about bounded-size objects to subsets of random sets. However, many difficult open problems remain regarding larger objects. It has long been known that many techniques, such as the regularity method, that are useful for embedding small graphs are also helpful for embedding large graphs of bounded degree. One might then expect that we can transfer such results to the sparse random setting in the same way that we have done for small graphs. This program has been partially successful (see, for example,~\cite{LS12} for work on a sparse random analogue of Dirac's theorem on Hamilton cycles and~\cite{ABHKP25} for a sparse random analogue of the powerful blow-up lemma), though many challenges remain. Here we will discuss these successes and their limitations through the lens of size-Ramsey numbers. 

Given a graph $H$ and a natural number $r$, the \emph{size-Ramsey number} $\hat{r}(H; r)$ is the smallest number of edges in an $(H,r)$-Ramsey graph. A foundational result of Beck \cite{B83}, which has since been extended in many ways (see, for example,~\cite{LPY25+}), says that for every integer $r \ge 2$ there exists a constant $C$ such that $\hat{r}(P_n;r) \leq C n$. Another key result in the area,  due to Kohayakawa, R\"odl, Schacht and Szemer\'edi~\cite{KRSS11}, says that if $H$ is any graph with $n$ vertices and maximum degree $\Delta$, then $\hat{r}(H;r) \leq n^{2 - \frac{1}{\Delta} + o(1)}$. That is, the size-Ramsey number of bounded-degree graphs is subquadratic in the number of vertices. More explicitly, they proved the following result about random graphs. 

\begin{thm} \label{thm:SizeRamsey}
For any integers $\Delta \geq 2$ and $r \geq 2$, there exists $C > 0$ such that if $p \geq C (\log N/N)^{1/\Delta}$, then the random graph $G_{N,p}$ with $N = C n$ a.a.s.~has the property that every $r$-colouring of the edges of $G_{n,p}$ has a colour class which contains every graph on $n$ vertices with maximum degree $\Delta$. 
\end{thm}

How tight is this bound? This contains two different questions, the first of which is: how small can we choose $p$ so that the random graph $G_{N,p}$ with $N = C n$ continues to be Ramsey with respect to every graph on $n$ vertices with maximum degree $\Delta$? Following discussions in~\cite{CNT22} and~\cite{AB25+}, we conjecture that Theorem~\ref{thm:SizeRamsey} continues to hold for all $p \geq C N^{-2/(\Delta+2)}$, which is the threshold coming from the R\"odl--Ruci\'nski theorem for $G_{N,p}$ to be Ramsey for $K_{\Delta + 1}$. That is, we conjecture that obtaining a monochromatic copy of $K_{\Delta + 1}$, the densest graph with maximum degree $\Delta$, is the main obstacle to obtaining monochromatic copies of all graphs with maximum degree $\Delta$.

\begin{con}
For any integers $\Delta \geq 2$ and $r \geq 2$, there exists $C > 0$ such that if $p \geq C N^{-2/(\Delta+2)}$, then the random graph $G_{N,p}$ with $N = C n$ a.a.s.~has the property that every $r$-colouring of the edges of $G_{n,p}$ has a colour class which contains every graph on $n$ vertices with maximum degree $\Delta$. 
\end{con}

This conjecture was confirmed for $\Delta = 3$ in~\cite{CNT22} and for $\Delta = 4$, up to a $o(1)$ factor in the exponent, in~\cite{AB25+}, with the latter paper also improving the bound for general $\Delta$ to $p \geq n^{-1/(\Delta - 1) + o(1)}$. However, even for $\Delta = 3$, this is not the end of the story, since one is not obliged to use vanilla random graphs when studying size-Ramsey numbers. By making use of a different random graph model with locally-dense spots, Dragani\'c and Petrova~\cite{DP25} were able to improve the upper bound on the size-Ramsey number of cubic graphs with $n$ vertices from the $Cn^{8/5}$ of~\cite{CNT22} to $n^{3/2 + o(1)}$. But this does not rule out the possibility that the bound could still be much lower. A first hurdle would be to show that there exists $\epsilon > 0$ such that the size-Ramsey number of cubic graphs with $n$ vertices is at most $C n^{3/2 - \epsilon}$, but the following key problem, first asked by R\"odl and Szemer\'edi~\cite{RSz00}, remains wide open. 

\begin{prob}
Prove or disprove that there is a constant $\epsilon > 0$ and an infinite sequence of cubic graphs $H$ with $\hat{r}(H;2) \ge v(H)^{1+\epsilon}$.
\end{prob}

At present, the best known lower bound for the size-Ramsey number of a cubic graph with $n$ vertices, due to Tikhomirov~\cite{T24} and building on earlier work of R\"odl and Szemer\'edi~\cite{RSz00}, stands at $n e^{c \sqrt{\log n}}$ for some $c > 0$.

\section{Combinatorial theorems relative to pseudorandom sets}
 
A pseudorandom subset of a set $X$ is one which behaves, in some sense, like a random subset. Usually, this is quantified through some family of statistics, such as near uniform distribution across some collection of subsets or containing asymptotically the same number of copies of some particular substructures as a random subset. Though we will not follow up on it here, one remarkable, and often extremely useful, phenomenon, explored in some depth by Chung, Graham and Wilson~\cite{CGW89} under the name of quasirandomness, is that many such pseudorandom properties are quantitatively related. 

For our purposes here, we will first be interested in a particular form of pseudorandomness for graphs, originating in work of Thomason~\cite{T87a, T87b}. 
We say that a graph $G$ is \emph{$(p, \beta)$-jumbled} if
\[|e(X, Y) - p|X||Y|| \leq \beta \sqrt{|X||Y|}\]
for any $X, Y \subseteq V(G)$. The random graph $G_{n,p}$ satisfies this condition with $\beta \le C \sqrt{pn}$, which is essentially best possible, though there are also many explicit examples which meet this bound. We refer the interested reader to the excellent survey~\cite{KS06} for more on this and on the basic properties of pseudorandom graphs, as well as useful descriptions of many important families of such graphs.

Given a graph property $\mathcal{P}$, an integer $n$ and a density $p$, who may ask for a determination of those values of $\beta$ such that any $(p, \beta)$-jumbled graph on $n$ vertices satisfies $\mathcal{P}$. For example, it is known that for any integer $t \geq 3$ there exists $c > 0$ such that if $\beta \leq c p^{t-1} n$, then any $(p, \beta)$-jumbled graph on $n$ vertices contains a copy of $K_t$. For $t = 3$, this condition is known to be tight, as shown by an example of Alon \cite{A94} (see also~\cite{C17, K11}), and it is a major open problem, with profound implications in Ramsey theory~\cite{MV24}, to show that it is also tight for all $t \ge 4$ (see~\cite{BIP20, MP22} for the best current bounds). Even a resolution of the following first open case would be of considerable interest.

\begin{prob}
Prove or disprove that there is a constant $C$ and an infinite sequence of $K_4$-free $(p, \beta)$-jumbled graphs on $n$ vertices with $p \ge C^{-1} n^{-1/5}$ and $\beta \leq C \sqrt{pn}$.
\end{prob}

For the problem of transferring combinatorial theorems to subgraphs of pseudorandom graphs, a general method, which they call densification, was developed by Conlon, Fox and Zhao in \cite{CFZ14a}. As an example of their results, we have the following pseudorandom analogue of Theorem~\ref{thm:removal-Gnp}.

\begin{thm} \label{thm:PseudoRemoval}
For any integer $t$ and any $\epsilon > 0$, there exist positive constants $\delta$ and $c$ such that if $\beta \leq c p^t n$, then any $(p, \beta)$-jumbled graph $G$ on $n$ vertices has the property that any subgraph of $G$ containing at most $\delta p^{\binom{t}{2}} n^t$ copies of $K_t$ may be made $K_t$-free by deleting at most $\epsilon p n^2$ edges. 
\end{thm}

This is in fact a special case of a more general statement that applies to all graphs and not just cliques (see also~\cite{ABSS20} for subsequent work with  improved quantification in certain cases). Moreover, with the same conditions on $\beta$, it is possible to prove analogues of Theorems~\ref{thm:RR} and~\ref{thm:randomTuran} and much more besides, in line with the transference theorems for random sets outlined in the previous section. However, it is worth stressing that the methods used for pseudorandom graphs are very different, and necessarily so, to those used in the random case, since, for instance, one does not have access to the union bound in the pseudorandom setting.  

There is still a gap in these results, even for triangles. For $t = 3$, Theorem~\ref{thm:PseudoRemoval} (which in this case was first proved by Kohayakawa, R\"odl, Schacht and Skokan \cite{KRSS10}) says that if $\beta \leq c p^3 n$, then the triangle removal lemma holds for subgraphs of a $(p, \beta)$-jumbled graph on $n$ vertices. However, it may well be that $\beta \leq c p^2 n$ is sufficient. That is, it may be that almost as soon as triangles are guaranteed to appear in pseudorandom graphs, they appear in such numbers that even the triangle removal lemma holds. 

\begin{prob} \label{prob:PseudoRemoval}
Prove or disprove that for any $\epsilon > 0$ there exist positive constants $\delta$ and $c$ such that if $\beta \leq c p^2 n$, then any $(p, \beta)$-jumbled graph $G$ on $n$ vertices has the property that any subgraph of $G$ containing at most $\delta p^3 n^3$ copies of $K_3$ may be made $K_3$-free by deleting at most $\epsilon p n^2$ edges. 
\end{prob}

It is known that $\beta \leq c p^2 n$ is a sufficient condition for other combinatorial properties of triangles, including the Ramsey property~\cite{CFZ14a} and the Tur\'an property~\cite{SSzV05}, to hold relative to $(p, \beta)$-jumbled graphs with $n$ vertices. However, a problem analogous to Problem~\ref{prob:PseudoRemoval} also remains open for the stability property, that is, for showing that triangle-free subgraphs of a pseudorandom graph $G$ with at least $(\frac12 - o(1)) e(G)$ edges must be within $o(e(G))$ edges of being bipartite. 
 
An extension of these results to hypergraphs was carried out in the follow-up paper~\cite{CFZ15}. For hypergraphs, we need to take a rather different viewpoint on pseudorandomness, as there is no clear analogue of the jumbled property (or the spectral theory from which it naturally arises via the celebrated expander mixing lemma~\cite{AC88}). Instead, if, given a $k$-uniform hypergraph $H$, we wish to prove a $H$-removal lemma relative to a much larger $k$-uniform hypergraph $G$, we need to assume that $G$ contains asymptotically the same number of copies as in a random hypergraph with the same density of each member of a certain family of small fixed subgraphs derived from $H$. Under such a condition, referred to as the $H$-linear forms condition, it is possible to prove a relative hypergraph removal lemma.  
Informally, this states that if $G$ satisfies the $H$-linear forms condition, then any subgraph of $G$ containing a $o(1)$-proportion of the copies of $H$ in $G$ can be made $H$-free by removing $o(e(G))$ edges. 
Since the precise statement is somewhat complicated, we refer the interested reader to~\cite{CFZ15} or the survey paper~\cite{CFZ14b} for more details.

One corollary of this result is an alternative proof of a key step in the proof of the Green--Tao theorem~\cite{GT08} that the primes contain arbitrarily long arithmetic progressions. Their proof can be viewed as having two main steps. The first step is to show that a certain set of almost primes, within which the primes form a dense subset, are pseudorandom in an appropriate sense. The second step is to show that Szemer\'edi's theorem continues to hold relative to such pseudorandom sets. In particular, since the primes are a dense subset of such a pseudorandom set, they must contain arbitrarily long arithmetic progressions.

In their paper, Green and Tao define a pseudorandom subset of the integers to be one that satisfies two different conditions, the linear forms condition and the correlation condition. Through the same mechanism that allows one to derive Szemer\'edi's theorem from the hypergraph removal lemma, the relative hypergraph removal lemma in~\cite{CFZ15} gives an alternative proof of their second step, the relative Szemer\'edi theorem, that only requires a linear forms condition. This simplification has been important to many of the subsequent results in the area, including, for instance, work of Tao and Ziegler~\cite{TZ18} on polynomial progressions in the primes.

It may be that even more can be derived from this approach. One particularly tantalising open problem is to show that there are infinitely many $3$-term arithmetic progressions of Friedlander--Iwaniec primes. A \emph{Friedlander--Iwaniec prime} is one of the form $x^2 + y^4$, so called because a seminal result of Friedlander and Iwaniec~\cite{FI98} says that, despite their rather low density, there are infinitely many such primes. In order to show that there are infinitely many $3$-term progressions of these primes using the approach taken above for ordinary primes, one would, besides the (likely considerable) additional difficulties on the number-theoretic side, need to show that a triangle-removal lemma holds down to density $p = n^{-1/4 - \epsilon}$ for some $\epsilon > 0$. At the moment, we know, from Theorem~\ref{thm:PseudoRemoval}, how to get down to roughly $n^{-1/5}$ and Problem~\ref{prob:PseudoRemoval} essentially asks whether it is possible to get down to roughly $n^{-1/3}$. Thus, settling that problem might be a first step towards the following one.

\begin{prob}
Prove that there are infinitely many $3$-term arithmetic progressions of Friedlander--Iwaniec primes.
\end{prob}
 
\section{Combinatorial theorems relative to extremal sets}
 
Given a natural number $n$ and a graph $H$, the \emph{extremal number} $\textrm{ex}(n, H)$ is the largest number of edges in an $H$-free graph with $n$ vertices. This is a much-studied function (see, for example,~\cite{FS13}), with the basic result, the Erd\H{o}s--Stone--Simonovits theorem that we already mentioned earlier, stating that
\[\textrm{ex}(n,H) = \left(1 - \frac{1}{\chi(H) - 1} + o(1)\right)\binom{n}{2}.\]
For $\chi(H) \ge 3$, this gives a fairly complete answer to the problem of determining $\textrm{ex}(n, H)$, but for $\chi(H) = 2$, that is, when $H$ is bipartite, it only shows that $\textrm{ex}(n, H) = o(n^2)$. 

A better bound was already given by K\H{o}v\'ari, S\'os and Tur\'an~\cite{KST54} in 1954, who showed that for any bipartite graph $H$ there exist $C, \delta > 0$ such that $\textrm{ex}(n, H) \le C n^{2-\delta}$. There are still relatively few bipartite graphs for which a more precise bound is known. However, the advent of the random algebraic method introduced by Bukh~\cite{B15} and developed further by Conlon~\cite{C19} and Bukh and Conlon~\cite{BC18} has added many additional cases. For example, the rational exponents conjecture of Erd\H{o}s and Simonovits~\cite{E81}, which states that for every rational number $r \in [1,2]$ there is a graph $H$ such that $\textrm{ex}(n, H) = \Theta(n^r)$, is now known for many exponents (see, for example,~\cite{CJ22, JQ23}) and a variant allowing for finite families instead of single graphs is completely resolved~\cite{BC18}.

The simplest well-understood case is when $H = C_4$, the cycle of length four. Here, it has been known~\cite{E38} since the 1930s that $\textrm{ex}(n, C_4) = \Theta(n^{3/2})$. Moreover, many of the constructions that give the lower bound are pseudorandom. Accordingly, in keeping with the theme of the last section, one might ask whether combinatorial theorems hold relative to $C_4$-free graphs. The answer is yes, at least to some extent. Conlon, Fox, Sudakov and Zhao~\cite{CFSZ21} studied this problem in some detail, proving several removal-type statements in this context. When it applies, a sparse $C_5$-removal lemma says that any subgraph of an $n$-vertex random or pseudorandom graph with density $p$ which contains $o(p^5 n^5)$ copies of $C_5$ can be made $C_5$-free by removing $o(p n^2)$ edges. Bearing in mind that the density of extremal $C_4$-free graphs is roughly $p = 1/\sqrt{n}$, the first main result of~\cite{CFSZ21} says that a $C_5$-removal lemma holds relative to $C_4$-free graphs at that $p$.

\begin{thm}
Every $n$-vertex $C_4$-free graph with $o(n^{5/2})$ copies of $C_5$ can be made $C_5$-free by removing $o(n^{3/2})$ edges. 
\end{thm}

Another of their removal lemmas has no explicit $C_4$-free assumption in the statement, though it is in a sense implicit, and the conclusion is about deleting copies of $K_3$ rather than $C_5$. 

\begin{thm}
Every $n$-vertex graph with no $C_5$ can be made $K_3$-free by deleting $o(n^{3/2})$ edges.
\end{thm}

An interesting corollary of this latter theorem concerns hypergraphs. To state it, we need to know that a \emph{(Berge) cycle} of length $r \ge 2$ in a hypergraph is an alternating sequence of distinct vertices and edges $v_1, e_1, \dots, v_r, e_r$ such that $v_i, v_{i+1} \in e_i$ for each $i$, where indices are taken modulo $r$. For example, a $2$-cycle consists of a pair of edges intersecting in a pair of distinct vertices. The \emph{girth} of a $k$-uniform hypergraph is then the length of the shortest cycle in the hypergraph.

\begin{thm} \label{thm:girth5}
Every $k$-uniform hypergraph on $n$ vertices of girth greater than $5$ has $o(n^{3/2})$ edges.
\end{thm}

We do not know if this is close to tight, but conjecture that it is.

\begin{con}
There are $3$-uniform hypergraphs of girth $6$ with $n$ vertices and $n^{3/2 - o(1)}$ edges.
\end{con}

Another corollary of these results (see also~\cite{P22} for an alternative approach) concerns finding solutions to linear equations within Sidon sets, which are in a sense the number-theoretic analogues of $C_4$-free graphs. Formally, a subset of $[n]$ is a \emph{Sidon set} if it contains no non-trivial solutions to the equation $x_1 + x_2 = x_3 + x_4$. It is known that such a Sidon set can have size at most $(1 + o(1))\sqrt{n}$ and also that there are Sidon sets of asymptotically this size (though determining the behaviour of the second-order term remains a deep and interesting open problem~\cite{BFR23}). 

\begin{thm} \label{thm:Sidon}
The maximum size of a Sidon subset of $[n]$ without a solution in distinct variables to the equation
\[x_1 + x_2 + x_3 + x_4 = 4x_5\] 
is at most $o(\sqrt{n})$ and at least $n^{1/2-o(1)}$.
\end{thm}

Here, we are simultaneously avoiding
\begin{itemize}
\item[(a)] non-trivial solutions to the Sidon equation $x_1 + x_2 = x_3 + x_4$ and
\item[(b)] distinct variable solutions to the linear equation $x_1 + x_2 + x_3 + x_4 = 4x_5$.
\end{itemize}
As mentioned above, there exist Sidon sets of size $(1+o(1))\sqrt{n}$ and, by a straightforward modification of Behrend’s construction~\cite{B46} of large sets without $3$-term arithmetic progressions, there are also sets of size $n^{1-o(1)}$ avoiding (b). What Theorem~\ref{thm:Sidon} says is that by simultaneously avoiding non-trivial solutions to both equations, the maximum size is substantially reduced.

In relation to extremal numbers, the well-known compactness conjecture of Erd\H{o}s and Simonovits~\cite{ES82} states that, for every finite family $\mathcal{F}$ of graphs, $\textrm{ex}(n, \mathcal{F}) \geq c_{\mathcal{F}} \min_{F \in \mathcal{F}} \textrm{ex}(n, F)$ for some constant $c_{\mathcal{F}} > 0$, where $\textrm{ex}(n, \mathcal{F})$ is the largest number of edges in an $n$-vertex graph not containing any graph from $\mathcal{F}$. The analogous statement has long been known to be false for $k$-uniform hypergraphs with $k \ge 3$, but, outside of some trivial counterexamples, remains open for graphs. Theorem~\ref{thm:Sidon} shows that it also fails for linear equations.
 
In a subsequent paper, Conlon, Fox, Sudakov and Zhao~\cite{CFSZ22} studied to what extent their techniques could be extended to prove counting (and thereby removal) lemmas for other $F$ besides $C_5$ relative to $C_4$-free graphs. Their result gives a recursive class of graphs, formed using ``islands'' and ``bridges'', that are countable in this sense. However, their result is likely far from complete and it remains an interesting open problem to determine exactly which graphs are countable. 
 
 \tikzstyle{v}=[circle, fill=black, inner sep = 1pt]
 
 \begin{prob} \label{open:dod-pet}
	 Which graphs are countable relative to $C_4$-free graphs? In particular, are the dodecahedral and Petersen graphs pictured below countable? 
	 
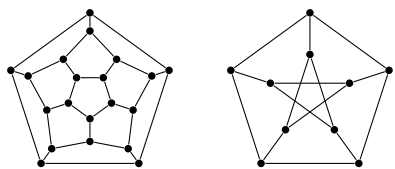
\[
	 \begin{tikzpicture}[scale=.2,baseline]
	 	\foreach \i in {0,1,2,3,4}{
	 		\node[v] (a\i) at (-90+72*\i:1.5) {};
	 		\node[v] (b\i) at (-90+72*\i:3) {};
	 		\node[v] (c\i) at (-54+72*\i:4.3) {};
	 		\node[v] (d\i) at (-54+72*\i:5.5) {};
	 	}
	 	\foreach \i[evaluate=\i as \j using {int(Mod(\i-1,5))}] in {0,1,2,3,4}{
	 		\draw (a\i) -- (a\j);
	 		\draw (a\i) -- (b\i);
	 		\draw (b\i) -- (c\i);
	 		\draw (b\i) -- (c\j);
	 		\draw (c\i) -- (d\i);
	 		\draw (d\i) -- (d\j);
	 	}
	 \end{tikzpicture}
	 \qquad 
 	 \begin{tikzpicture}[scale=.55,baseline]
	 	\foreach \i in {0,1,2,3,4}{
 		\node[v] (a\i) at (90+72*\i:1) {};
	 		\node[v] (b\i) at (90+72*\i:2) {};
	 	}
	 	\foreach \i[evaluate=\i as \ia using {int(Mod(\i+2,5))},evaluate=\i as \ib using {int(Mod(\i+1,5))}] in {0,1,2,3,4}{
	 		\draw (a\i) -- (b\i);
	 		\draw (a\i) -- (a\ia);
	 		\draw (b\i) -- (b\ib);
	 	}
	 \end{tikzpicture}
	 \]
 \end{prob}
 
Another, arguably more important, open problem is to extend the techniques of~\cite{CFSZ21, CFSZ22} to $C_{2k}$-free graphs for $k \ge 3$. In particular, when $k = 3$, one can ask whether a $C_7$-removal lemma holds relative to $C_6$-free graphs, where here the natural density is at $p = n^{-2/3}$.  
 
\begin{prob}
Prove that every $n$-vertex $C_6$-free graph with $o(n^{7/3})$ copies of $C_7$ can be made $C_7$-free by removing $o(n^{4/3})$ edges.
\end{prob}

If such a statement could be proved, it is likely that one could also prove an analogue of Theorem~\ref{thm:girth5} saying that $k$-uniform hypergraphs on $n$ vertices of girth greater than $7$ have $o(n^{4/3})$ edges and an analogue of Theorem~\ref{thm:Sidon} regarding solutions of certain linear equations relative to $B_3$-sets.

\address{Department of Mathematics, California Institute of Technology\\
 Pasadena, CA 91125, USA.\\
\email{dconlon@caltech.edu}}

\end{document}